# Harmonic Analysis

Vladimir Clue
*Vladimir Clue*
vclue@earthlink.net

*I thank Natasha Us, William Kelleher and Anatoly Levin for helping me to make out this paper.*

## Abstract

*This paper describes a method of calculating the transforms, currently obtained via Fourier and reverse Fourier transforms. The method allows calculating efficiently the transforms of a signal having an arbitrary dimension of the digital representation by reducing the transform to a vector-to-circulant matrix multiplying.*

*There is a connection between harmonic equations in rectangular and polar coordinate systems. The connection established here and used to create a very robust iterative algorithm for a conformal mapping calculation.*

*There is also suggested a new ratio (and an efficient way of computing it) of two oscillative signals.*

## 1. Introduction

The traditional way to transform a signal **U** i.e., to change amplitude-frequency characteristic (AFC) of **U** in some desired manner, is: 1) use Fourier transform to obtain the AFC of the signal; 2) apply some mechanism to recalculate the row coefficients; and 3) use inverse Fourier transform to obtain the transformed signal. Also, the traditional way to calculate values of a function **U** + i***V**, analytical inside a unit circle, from given values on the unit circle, so that new values resemble values of that function on some other concentric circle, is: 1) use Fourier transform to obtain row coefficients; 2) recalculate the coefficients to reflect that the new values are on the circle with radius 'r' and starting argument 'ψ'; and 3) use reverse Fourier transform to obtain the desired values.

The traditional way to transform a signal also has a well-known obstacle. The obstacle is the increasing relative cost of calculations that occurs with an increasing number of sampling points **N**. If **N** is a power of 2 the cost is **N*log$_2$ (N)**, but in other cases the cost is higher. A worst-case scenario is the absence of small prime dividers, which results in costs proportional to **N^2**.

I intend to prove here that one circulant matrix operator can do both transforms: transform the AFC, and transform from analytical function values on the unit circle to the values on the concentric circle with a different radius and starting point argument. In both transforms, only **U** is used as an input vector, while the result of the transform is a complex vector.

A very efficient algorithm exists for multiplying circulant matrixes by vectors without having to worry about a dimension of the matrix being either a power of 2 or having other small prime dividers. [1] Therefore, a direct transform with circulant matrix can be done more efficiently than via the Fourier's transform followed by the reverse Fourier transform.

## 2. Transform the AFC by Circulant Matrix

It is well known that the AFC of a digital signal **U** with **N** number of sampling points is

$$a_k = \frac{2}{N} * \sum_{t=0}^{N-1}(U_t \cos(2\pi \frac{kt}{N})),$$
$$b_k = \frac{2}{N} * \sum_{t=0}^{N-1}(U_t \sin(2\pi \frac{kt}{N})) \qquad (1)$$

It is also well known that the row coefficients for a function $W = U + iV$, analytical inside a unit circle are:

$$C_j = \frac{1}{N} * \sum_{t=0}^{N-1} W(\exp(2\pi i \frac{t}{N})) \exp(-2\pi i \frac{jt}{N})$$

**Corollary I**
If a circulant matrix has elements defined as:

$$M_{l,j}(\lambda_t, r, \psi) = \frac{2}{N} * \sum_{k=1}^{(N-1)/2}(\lambda_k r^k \exp(i(2\pi \frac{l-j}{N} + \psi)k) \qquad (3)$$

where $\lambda_t$, **r** and $\psi$ are real numbers (for convergence reason parameter 'r' has to be less than the radius of convergence of the row representing the analytical function[1]), then the result of the transform with the given matrix can be presented in two ways:

i)
$$(M*U)_j = \sum_{k=1}^{K}(a\_new_k \cos(2\pi \tfrac{kj}{N}) + b\_new_k \sin(2\pi \tfrac{kj}{N})) +$$
$$+ i*\sum_{k=1}^{K}(-b\_new_k \cos(2\pi \tfrac{kj}{N}) + a\_new_k \sin(2\pi \tfrac{kj}{N}))$$

where the new AFC relates to the original one as:
$$a\_new_k = \lambda_k r^k (a_k \cos(k\psi) + b_k \sin(k\psi))$$
$$b\_new_k = \lambda_k r^k (b_k \cos(k\psi) - a_k \sin(k\psi))$$

or

ii) $$(M*U)_j = \sum_{t=1}^{N-1} c\_new_t \exp(2\pi i \tfrac{jt}{N})$$

where the new coefficients relate to the original ones as:
$$c\_new_j = \lambda_j r^j c_j \exp(i\psi j)$$

It is worth noting that the transform "loses" a constant – zero's harmonic and also it "loses" $N/2$'s harmonic in the case where N is an even number. This transform always produces a result, as if the amplitude(s) of the harmonic(s) were zero.

**Proof:**
To prove Corollary I, I simply construct the operator **M** from two parts.
The first part is similar to the way in which a Szego kernel is constructed:

$$M1 = \left(\cos(k(2\pi \tfrac{l}{N}+\psi)), \sin(k(2\pi \tfrac{l}{N}+\psi))\right)*$$
$$*(\mu_k)\otimes\left(\cos(k(2\pi \tfrac{j}{N})), \sin(k(2\pi \tfrac{j}{N}))\right)^T$$

Here $\otimes$ means multiplying components so that corresponding cosine-s and sine-s are multiplied by corresponding coefficients $\mu_k$.
In the case where $\psi=0$, $\mu_k \equiv 2/N$, and k=1,..(N-1)/2, this part becomes the Szego kernel [2] except that it does not have a reconstructive property for zero's harmonic and N/2's harmonic, if N is an even number.
The second part is harmonically conjugated to the first part:
$$M2 = \left(\sin(k(2\pi \tfrac{l}{N}+\psi)), -\cos(k(2\pi \tfrac{l}{N}+\psi))\right)*$$
$$*(\mu_k)\otimes\left(\cos(k(2\pi \tfrac{j}{N})), \sin(k(2\pi \tfrac{j}{N}))\right)^T$$

---

[1] If the condition is satisfied, then the accuracy improves with increasing the number of sampling points, otherwise the result is unpredictable.

---

The elements of M1 and M2 are:

$$M1_{l,j} = \sum_{k=1}^{(N-1)/2} \mu_k (\cos(k(2\pi \tfrac{l}{N}+\psi))\cos(k(2\pi \tfrac{j}{N})) +$$
$$+\sin(k(2\pi \tfrac{l}{N}+\psi))\sin(k(2\pi \tfrac{j}{N})))$$

$$M2_{l,j} = \sum_{k=1}^{(N-1)/2} \mu_k (\sin(k(2\pi \tfrac{l}{N}+\psi))\cos(k(2\pi \tfrac{j}{N})) -$$
$$-\cos(k(2\pi \tfrac{l}{N}+\psi))\sin(k(2\pi \tfrac{j}{N})))$$

The elements of $(M1+i*M2)_{l,j}$, after applying trigonometric formulas for cosine-s and sine-s of angle differences, are:

$$(M1+i*M2)_{l,j} = \sum_{k=1}^{(N-1)/2} \mu_k [\cos(k(2\pi \tfrac{l}{N}+\psi) - k(2\pi \tfrac{j}{N})) +$$
$$+ i\sin(k(2\pi \tfrac{l}{N}+\psi) - k(2\pi \tfrac{j}{N}))]$$

Now, if one applies Euler's formula for complex exponent and substitutes $\mu_k = \tfrac{2}{N}\lambda_k r^k$, one will arrive at (3). This discussion proves the Corollary I for representation given in (1). Proof of the corollary I for representation (2) is very similar, once one remembers that V is harmonically conjugated to U. [3]

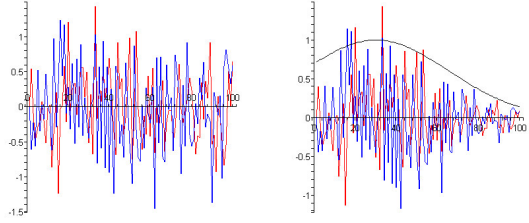

Fig. 1
Fig.1 shows the AFC of a signal and the AFC of the signal transformed on Maplesoft platform by matrix

M($\lambda$, 1, 0),

Where:
$\lambda_k =$ `exp(-(0.3*K-k)^2/K^2*4)`,
`K = (Dimension(signal)-1)/2, and`
`k = 1..K.`

There are many cases when the polynomial expression (3) can be factored. While some of the cases can be applied to calculate customized filters efficiently, the one particular case allows evaluating an analytical function given by its real part values in the unit circle.

## 1.1 Evaluating Analytical Functions.

The operator given by (3) is very useful for evaluating analytical functions. Indeed, if one sets $\lambda_k \equiv 1$ and uses the formula for a sum of geometrical progression, then the operator will become

$$MA_{l,j}(r,\psi) = \tfrac{2}{N} r \exp(2\pi i \tfrac{l-j}{N} + i\psi) *$$
$$* r^{(N-1)/2} \frac{\exp((2\pi i \tfrac{l-j}{N} + i\psi)(N-1)/2) - 1}{\exp(2\pi i \tfrac{l-j}{N} + i\psi) - 1} \quad (4)$$

which allows the evaluation of an analytical function based upon its real part values on the same circle. New values of an analytical function on the circle having r as a ratio to the initial radius and new starting argument having $\psi$ as a difference from the initial starting point can be obtained by applying the "MA" operator on the real part of the analytical function.

## 1.2 Obtaining values of harmonically conjugated function.

Of special interest is the case when one needs to obtain values of a harmonically conjugated function on a circle by given real part values on the circle. I simplified the expression by setting values r=1 and $\psi$=0 and extracting the imaginary part. After fairly simple steps, I obtained:

$$MH_{l,j} = \frac{\cos(\tfrac{\pi}{N}(l-j)) - \cos(\tfrac{\pi}{N}(l-j))^{1-\text{mod}(N,2)} * (-1)^{\text{mod}(l-j,2)}}{N \sin(\tfrac{\pi}{N}(l-j))}$$

where the main diagonal elements apparently are nulls.

## 2. Connection between harmonic equations in rectangular and polar coordinate systems

It is well known that if $W(z) = U(x,y) + iV(x,y)$ is an analytical inside of the unit circle one-listed function, then $U(x,y)$ and $V(x,y)$ are harmonic functions and they are harmonically conjugated. Therefore, the values of those functions on the unit circle are connected as $V = MH(U)$. Later on I will denote "MH" operator as "~" – operator. It is also well known that arg(z) is harmonically conjugated to ln(|z|). Corollary II states the more generalized fact.

**Corollary II**

If $W(z) = U(x,y) + iV(x,y)$ is an analytical inside of the unit circle one-listed function, then the radius vector's norm and argument of the image of the unit circle produced by the given analytical function are connected as

$$\varphi(\alpha) - \alpha = \sim (Ln[\sqrt{norm(W(\alpha) - W(0))}]) \quad (5)$$

where $\alpha = \arg(z)$ and $\varphi = \arg(W)$.
Furthermore, the tangent angle $\phi(\alpha)$ to boundary curve in $W(\alpha)$ and modulus of derivative are connected as

$$\phi(\alpha) - (\alpha + \tfrac{\pi}{2}) = \sim (Ln[\sqrt{norm(W(\alpha)')}]) \quad (6)$$

**Proof:**
Row representation of the analytical, one-to-one function W(z) is

$$W(z) = \sum_{n=0}^{\infty} C_n z^n \quad (7)$$

Function $W1 = \frac{W(z) - W(0)}{z}$ has the following row representation:

$$W1(z) = \sum_{n=0}^{\infty} C_{n+1} z^n \quad (8)$$

If row (7) is absolutely convergent, then row (8) is also absolutely convergent inside radius of convergence of (7) since it has asymptotically smaller coefficients for same powers of z; therefore the function W1 is also analytical. W1 has no zeros: in z=0, it is equal to the derivative of W and if it was equal to zero in z≠0, then (z*W1) would have two zero values inside the radius of convergence, which would make it not one-to-one; therefore, Ln(W1) is analytical and its real and imaginary parts are harmonically conjugated. If one expresses W1 in polar coordinates and take the natural logarithm, one will arrive at expression (5).
Expression (6) is obtained similarly. The logarithm of a derivative expressed in polar coordinates is an analytical function, so its real and imaginary parts are harmonically conjugated, and (6) states just that.

### 2.1 Conformal Mapping

Equations (5) and (6) are powerful tools for solving Riemann's task of finding conformal mapping from unit disk to simple-connected area surrounded by Jordan's curve. [4] Indeed, since ~(~(x)) = -x, then (6) implies that

$$\sqrt{norm(W(\alpha)')} = \exp(\sim(-(\phi(\alpha) - (\alpha + \tfrac{\pi}{2}))))$$

So, given an initial approximation, the density distribution can be corrected. The best results can be obtained by combining both (5) and (6). Using (6) is especially important, since (5) does not allow control over W(0), which in (6) one can set according to one's needs.

The combining can be achieved by calculating dissatisfaction function values, which provide a ratio by which the density distribution has to be changed. The task is to find reparametrisation t(s) iteratively adjusting it according to the dissatisfaction function, to eventually satisfy equations (5) and (6).

In the first step, t(s) is taken as t(s) = s. The iteration is done via consecutive steps of calculating the non-satisfaction function in the following way:

```
For j ∈ [0; N-1]
  ns(j) = exp [~((a_(j)-(π/2+j*2π/N))    \
                 -(a(j)-j*2π/N))]*       \
          sqrt[norm(w'(j))/norm(w (j)-w0)];
ns *= (N/sqrt(norm (ns)),
```

where a(j) = arg(w(j)-w0), a_(j)-j*2π/N = arg(w'(j)). When the Ln(ns(j)) is less than a predetermined tolerance, the points with coordinates {u(t(s(j))), v(t(s(j)))} resemble an image of the base grid that is reflected by a function that is analytical (having the predetermined tolerance) inside the unit circle.

This algorithm converges exponentially. Typically, only several iterations are needed to obtain a solution with a very high accuracy.

The solution of Riemann's task therefore, can be obtained at a time proportional to N*Ln(N) number of operations, although the evaluation on the N*N mesh will require N*N*Ln(N) number of operations.

Pictured below are some practical results produced by the new method implemented in Maple 9.5 and VC++. The number of iterations in all cases was no more then 10. It took 4 seconds for a PC to compute and draw the last one (Fig. 5), calculated for extra high accuracy with 2000 sampling points.

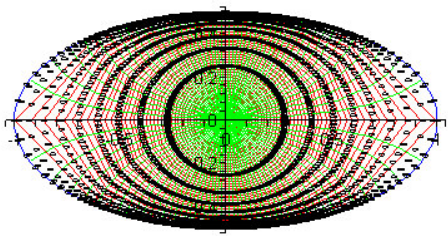

Fig. 2

The green-red maps in Fig. 2 and Fig. 3 are the result of Maple's 'conformal' function, black-dotted is the result of transform given by matrix (4). On the right side of conformal images in Fig. 4 and Fig. 5 are the plots of the reparametrisation and non-satisfaction functions.

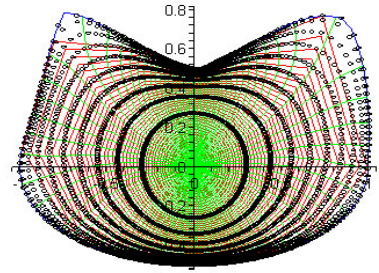

Fig. 3

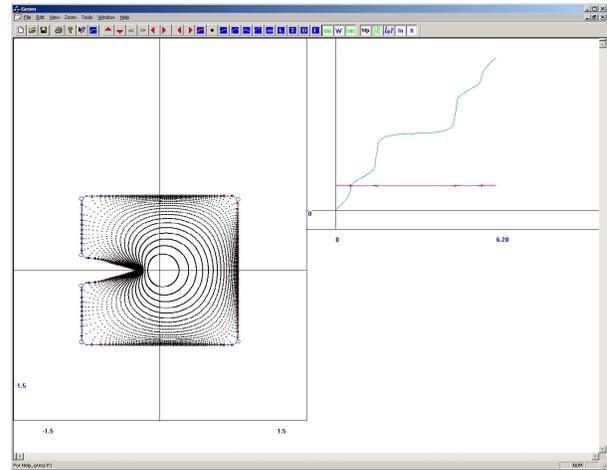

Fig. 4

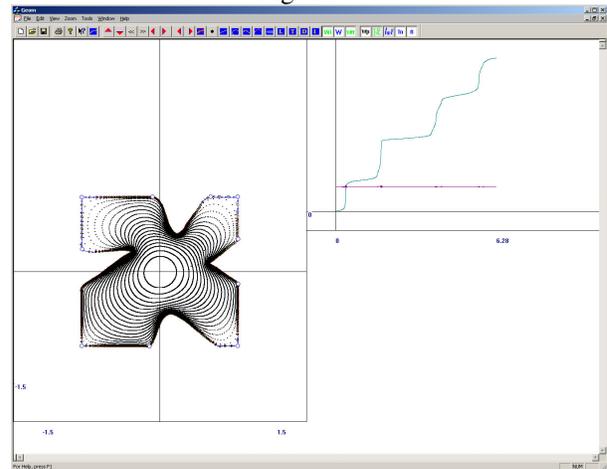

Fig. 5

Investigation of the convergence of the method is not yet complete; Therefore, I limit myself here to just an illustration of a very basic idea of why it works: Let -δ*Sin(kt+ψ) be a deviation component of frequency k, that is if ideal solution is: $|f'| = \exp(\sim(\alpha - \varphi))$, then for the deviated distribution the ratio is $\exp(\sim(\alpha-\varphi_1))/|f'_1|$, where $\varphi_1 = \varphi - \delta\sin(k\alpha+\psi))$. The ratio gives correcting multiplier: $ns(\alpha) = \exp(-\delta\cos(k\alpha+\psi))*|f'|/|f'_1|$

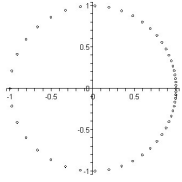

Fig. 6

Fig. 6 illustrates the case when k = 1 and ψ = 0.

## 2.2 Reinstating Wave Function

Expression (5) can also be used to reinstate wave function. Only the norm of the wave function is available to us as the result of experiments in quantum mechanics – it is the probability density distribution. [5] From the norm known on the circle, flat wave function, according to (5) is

$$\Psi(\xi) = \sqrt{|\Psi(\xi)|^2} * \exp(i \sim [Ln(\sqrt{|\Psi(\xi)|^2})])$$

where Ψ(ζ) is the wave function, and $(|\Psi(\zeta)|)^2$ denotes its norm.

## 3. Harmonic Covariation, Correlation

### 3.1 Harmonic Covariation

For two oscillative function U(t) and U(t) integrated on interval [0; 2π], I define harmonic covariation as:

$$U \sim V = \frac{1}{\pi}(\int_0^{2\pi}(U-mean(U))(V-mean(V))dt + \\ + i\int_0^{2\pi}(U-mean(U))[\sim(V-mean(V))]dt)$$

The properties of the tilde operator are:

1) (α + i*β) ~ u = α*u + β*(~u),
2) λ ~ u = u ~ conj(λ),
3) (λ + μ)~u = λ ~ u + μ ~ u,
4) u ~ v = conj(v ~ u), and
5) u~(λ~x + μ~y) = conj(λ)*(u~x) + conj(μ)*(u~y),

where α, β are real numbers; λ and μ are complex numbers; and u, v, x and y are functions defined on interval [0; 2π]

Property number 2 is not obvious; an illustrative example is as follows:

[ sin(a) ~ cos(a) ] ~ sin(a) = i ~ sin(a)  &

sin(a) ~ [cos(a) ~ sin(a)] = sin(a) ~ (-i) ⇒
i~sin(a) = sin(a)~(-i)

The essence of tilde-multiplying a complex number to a vector is described by:

$(a+ib)*U = aU + b*(\sim U)$

So, tilde-multiplying complex number to real-valued vector is a real-valued vector.

### 3.2 Harmonic Correlation

I define Harmonic Correlation HC(U,V) as:

$$HC(U,V) = \frac{U \sim V}{sqrt(norm(U-mean(U))norm(V-mean(V)))}$$

This ratio may be illustrated as follows. The HC of W(z), analytical inside unit circle, and W(z), multiplied by a constant equal to exp(iφ), for some values of φ are:

HC(Real(W(α)), Real(W(α)exp(0))) = 1
HC(Real(W(α)), Real(W(α)exp($i\frac{\pi}{2}$))) = i
HC(Real(W(α)), Real(W(α)exp(iπ))) = −1
HC(Real(W(α)), Real(W(α)exp($i\frac{3\pi}{2}$))) = −i

Therefore, harmonic correlation is useful for calculating the ratio of signals in cases when only one of the components of a complex valued variable is available for observation and/or measurement. One example could be the electric volume of the complex electric impedance.

In another example, the method of harmonic correlation was used to sort market data to find companies with share price behavior that is most similar to some particular company's share price behavior.

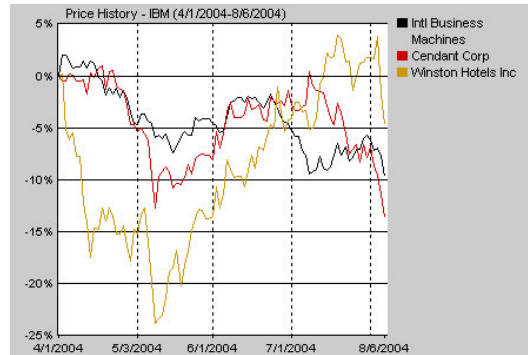

Fig. 7

In Fig. 7 and Fig. 8, market data is shown that was sorted by using the standard correlation coefficient (Fig. 7) and harmonic correlation coefficient (Fig. 8) to find the stocks whose price history most closely tracked that of International Business Machines Corporation (IBM).

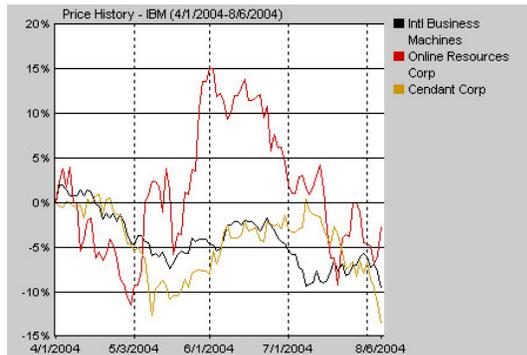

Fig. 8

The data of fifteen-minute intervals for a two-week time period was used for computing the ratios. The two closest tracking companies are shown along with IBM. The smaller the variation and the greater similarity in behavior in the share price data among those stocks selected using the method of Harmonic Correlation demonstrate the advantage of the Harmonic Correlation as a ratio between oscillative signals.